

\baselineskip=14pt
\parskip=10pt
\def\halmos{\hbox{\vrule height0.15cm width0.01cm\vbox{\hrule height
  0.01cm width0.2cm \vskip0.15cm \hrule height 0.01cm width0.2cm}\vrule
  height0.15cm width 0.01cm}}

\magnification=\magstephalf

\def\1{{\overline{1}}}
\def\2{{\overline{2}}}
\parindent=0pt
\overfullrule=0in

\def\frac#1#2{{#1 \over #2}}
\centerline
{\bf  Surprising Relations Between Sums-Of-Squares of Characters of the Symmetric Group }
\centerline
{\bf Over Two-Rowed Shapes and Over Hook Shapes }
\bigskip
\centerline
{\it By Amitai REGEV and Doron ZEILBERGER}

{\bf Abstract}: In a recent article, we noted (and proved) that the sum of the squares of the characters of the symmetric group,
$\chi^{\lambda}(\mu)$,
over all shapes $\lambda$ with two rows and $n$ cells and $\mu=31^{n-3}$, equals, surprisingly, to
1/2 of that sum-of-squares taken over all hook shapes  with $n+2$ cells and with $\mu= 321^{n-3}$.
In the present note, we show that this is only the tip of a huge iceberg! We will prove that if
$\mu$ consists of odd parts and (a possibly empty) string of {\it consecutive} powers of $2$, namely $2,4, \dots, 2^{t-1}$ for $t \geq 1$, then
the the sum of  $\chi^{\lambda}(\mu)^{2}$ over all two-rowed shapes $\lambda$  with $n$ cells, equals exactly $\frac{1}{2}$ times
the analogous sum of $\chi^{\lambda}(\mu')^{2}$ over all shapes $\lambda$ of {\it hook shape} with $n+2$ cells, and
where $\mu'$ is the partition obtained from $\mu$ by retaining all odd parts, but replacing the string
$2,4, \dots, 2^{t-1}$ by $2^t$.

Recall that the {\it Constant Term} of a {\it Laurent polynomial} in $(x_1, \dots, x_m)$ is the free term, i.e.
the coefficient of $x_1^0 \cdots x_m^0$. For example

$$
CT_{x_1,x_2} (x_1^{-3}x_2+x_1x_2^{-2}+5)=  5 \quad  .
$$

Recall that a {\it partition} (alias {\it shape}) of an integer $n$, with $m$ {\it parts} (alias {\it rows}),
is a non-increasing sequence of positive integers
$$
\lambda= (\lambda_1, \dots, \lambda_m) \quad,
$$
where $\lambda_1 \geq \lambda_2 \geq \dots \geq \lambda_m >0$, and $\lambda_1 + \dots + \lambda_m=n$.

If $\lambda=(\lambda_1, \dots, \lambda_m)$ and $\mu=(\mu_1, \dots, \mu_r)$ are partitions of $n$ with $m$ and $r$ parts,
respectively, then it easily follows from (7.8) (p. 114) in [M], that
the {\it characters}, $\chi^{\lambda}(\mu)$,  of the {\it symmetric group}, $S_n$, may be obtained
via the {\it constant term} expression
$$
\chi^{\lambda}(\mu) \, = \,
CT_{x_1, \dots , x_m} \,
\frac
{\prod_{1 \leq i < j \leq m} (1-\frac{x_j}{x_i}) \prod_{j=1}^r \left ( \sum_{i=1}^{m} x_i^{\mu_j} 
\right )}
{\prod_{i=1}^{m} x_i^{\lambda_i}} 
  \quad .
\eqno(Chi)
$$

As usual, for any partition $\mu$, $|\mu|$ denotes the sum of its parts, in other words, the integer that is
being partitioned.

In [RRZ] we considered  two quantities. Let $\mu_0$ be any partition with smallest part $\geq 2$.
The first quantity, that we will call henceforth
$A(\mu_0)(n)$, is the following  sum-of-squares over two-rowed shapes $\lambda$:
$$
A(\mu_0)(n):= \, \sum_{j=0}^{\lfloor n/2 \rfloor} \chi^{(n-j,j)}( \mu_0 1^{n-|\mu_0|})^2 \quad .
$$
[Note that in [RRZ] this quantity was denoted by $\psi^{(2)}(\mu_0 1^{n-|\mu_0|})$. ]

The second quantity was the sum-of-squares over {\it hook-shapes}
$$
B(\mu_0)(n):= \, \sum_{j=1}^{n} \chi^{(j,1^{n-j})}( \mu_0 1^{n-|\mu_0|} )^2 \quad.
$$
[Note that in [RRZ] this quantity was denoted by $\phi^{(2)}(\mu_0 1^{n-|\mu_0|})$.]

In [RRZ]  we developed algorithms for discovering (and then proving) 
closed-form expressions for these quantities, for any {\it given} (specific)
finite partition $\mu_0$ with
smallest part larger than one. In fact we proved that each such expression is {\it always} a multiple of
${{2n} \choose {n}}$ by a certain rational function of $n$ that depends on $\mu_0$.

Unless $\mu_0$ is very small, these rational functions turn out to be very complicated, but, inspired by
the OEIS([S]),  Alon Regev noted (and then it was proved in [RRZ]) the {\it remarkable} identity
$$
A(3)(n)=\frac{1}{2} B(3,2)(n+2) \quad .
$$

This lead to the following natural question:

Are there other partitions, $\mu_0$, such that there exists a partition, $\mu'_0$ with
$|\mu'_0|=|\mu_0|+2$, such that the ratio $A(\mu_0)(n)/B(\mu'_0)(n+2)$ is a constant? 

This lead us to write a new procedure in the Maple package

{\tt http://www.math.rutgers.edu/\~{}zeilberg/tokhniot/Sn.txt} \quad , \quad that accompanies [RRZ],

called  {\tt SeferNisim(K,N0)}, that searched for such pairs $[\mu_0,\mu'_0$]. 
We then used our {\it human}  ability for {\it pattern recognition} to notice that  all the successful pairs 
(we went up to $|\mu_0| \leq 20$) turned out to be such that  $\mu_0$ either consisted of only odd parts, and then $\mu'_0$ was
$\mu_0$ with $2$ appended, or, more generally $\mu_0$ consisted of odd parts together with a string of {\it consecutive}
powers of $2$ (starting with $2$), and $\mu'_0$ was obtained from $\mu_0$ by retaining all the odd parts but replacing the
string of powers of $2$ by a single power of $2$, one higher then the highest in $\mu_0$. In symbols,
we conjectured, (and later proved [see below], {\it alas}, by purely human means) the following:

{\bf Theorem:} Let $\mu_0$ be any partition of the form
$$
\mu_0 =Sort([a_1, \dots, a_s, 2,2^2, \dots, 2^{t-1}]) \quad,
$$
where
$$
a_1 \geq a_2 \geq \dots \geq a_s \geq 3 \quad,
$$
are all {\bf odd}, and $t \geq 1$ (if $t=1$ then $\mu_0$ only consists of odd parts). Define
$$
\mu'_0 =Sort([a_1, \dots, a_s, 2^t]) \quad.
$$
Then, for every $n \geq |\mu_0|$, we have
$$
A(\mu_0)(n) = \frac{1}{2} B(\mu'_0)(n+2) \quad .
$$
(For any sequence of integers, S, Sort(S) denotes that sequence sorted in non-increasing order.)

In order to prove our theorem we need to first recall, from [RRZ], the following {\bf constant-term} expression for
$B(\mu_0)(n)$.

{\bf Lemma 1}: Let $\mu_0=(a_1, \dots, a_r)$
$$
B(\mu_0)(n)= Coeff_{x^0}
\left [ 
\frac{(1+x)^{2n-2 -2(a_1 + \dots + a_r) }}{x^{n-1}} \cdot \prod_{i=1}^{r} (x^{a_i} -(-1)^{a_i}) (1 - (-1)^{a_i}x^{a_i} ) 
\right ] \quad .
$$

We need an analogous constant-term expression  for $A(\mu_0)(n)$. To that end, let's first
spell-out Equation $(Chi)$ for the two-rowed case, $m=2$, so that we can write $\lambda=(n-j,j)$.
We have, writing $\mu_0=(a_1, \dots, a_r)$,
$$
\chi^{(n-j,j)}(\mu_0 1^{n-|\mu_0|}) \, = \,
CT_{x_1, x_2} \,
\frac
{ (1-\frac{x_2}{x_1})  (x_1 + x_2)^{n-a_1- \dots -a_r}  \prod_{i=1}^r \left ( x_1^{a_i} + x_2^{a_i} \right )}
{x_1^{n-j} x_2^{j}}   \quad .
\eqno(Chi2)
$$
This can be rewritten as
$$
\chi^{(n-j,j)}(\mu_0 1^{n-|\mu_0|}) \, = \,
CT_{x_1, x_2} \,
\frac
{ (1-\frac{x_2}{x_1})  (1 + \frac{x_2}{x_1})^{n-a_1- \dots -a_r}  \prod_{i=1}^r \left ( 1+ (\frac{x_2}{x_1})^{a_j} \right )}
{(\frac{x_2}{x_1})^j }   \quad .
\eqno(Chi2')
$$
Since the {\it constant-termand} is of the form $P(\frac{x_2}{x_1})/(\frac{x_2}{x_1})^j$ for some {\it single-variable} polynomial $P(x)$, the
above can be rewritten, as
$$
\chi^{(n-j,j)}(\mu_0 1^{n-|\mu_0|}) \, = \,
Coeff_{x^0} \,
\frac
{ (1-x)  (1 + x)^{n-a_1- \dots -a_r}  \prod_{i=1}^r \left ( 1+ x^{a_i} \right )}
{x^j }   \quad .
\eqno(Chi2'')
$$

Note that the left side is {\it utter nonsense} if $j > \frac{n}{2}$, but the right side makes perfect sense. It is easy to see that
defining $\chi^{(n-j,j)}(\mu_0 1^{n-|\mu_0|})$ by the right side for  $j > \frac{n}{2}$, we get
$$
\chi^{(n-j,j)}(\mu_0 1^{n-|\mu_0|}) \, = \, - \chi^{(j,n-j)}(\mu_0 1^{n-|\mu_0|}) \quad .
$$ 

Let's denote the numerator of the constant-termand of $(Chi'')$,
namely 
$$
(1-x)  (1 + x)^{n-a_1- \dots -a_r}  \prod_{i=1}^r \left ( 1+ x^{a_i} \right ) \quad ,
$$
by $P(x)$, then
equation $(Chi2'')$ can be also rewritten as a {\it generating function}.
$$
P(x)=\sum_{j=0}^{n} \chi^{(n-j,j)}(\mu_0 1^{n-|\mu_0|}) \, x^j \quad .
$$
Since for any polynomial of a single variable, $P(x)=\sum_{j=0}^{n} c_j x^j$, we have
$$
\sum_{j=0}^{n} c_j^2 = Coeff_{x^0} \left [ P(x) P(x^{-1}) \right ] \quad ,
$$
we get
$$
\sum_{j=0}^{n}  \chi^{(n-j,j)}(\mu_0 1^{n-|\mu_0|})^2 =
$$
$$
Coeff_{x^0} 
\left [ \left ( (1-x)  (1 + x)^{n-a_1- \dots -a_r}  \prod_{j=1}^r \left ( 1+ x^{a_j} \right ) \right ) \cdot
\left ( (1-x^{-1})  (1 + x^{-1})^{n-a_1- \dots -a_r}  \prod_{j=1}^r \left ( 1+ x^{-a_j} \right ) \right )
\right ] \quad .
$$
$$
=
-Coeff_{x^0}
\left [
\frac{ (1-x)^2  (1 + x)^{ 2(n-a_1- \dots -a_r)}  \prod_{j=1}^r \left ( 1+ x^{a_j} \right )^2 }
{x^{n+1}} \, \right ] \quad .
$$
But since, by symmetry,
$$
\sum_{j=0}^{\lfloor \frac{n}{2} \rfloor }  \chi^{(n-j,j)}(\mu_0 1^{n-|\mu_0|})^2 = 
\frac{1}{2} \sum_{j=0}^{n}  \chi^{(n-j,j)}(\mu_0 1^{n-|\mu_0|})^2  \quad,
$$
we have

{\bf Lemma 2}: Let $\mu_0=(a_1, \dots, a_r)$ be a partition with smallest part larger than one, then
$$
A(\mu_0)(n)=
-\frac{1}{2} \, Coeff_{x^0} \left [
\frac{(1-x)^2  (1 + x)^{2(n-a_1- \dots -a_r)}  \prod_{j=1}^r \left ( 1+ x^{a_j} \right )^2 }
{x^{n+1}} \right ] \quad .
$$

We are now ready to prove the theorem. If $\mu_0=Sort (a_1, \dots, a_r, 2, \dots , 2^{t-1})$ then
$$
A(\mu_0)(n)=
- \frac{1}{2} \, Coeff_{x^0} \left [
\frac{(1-x)^2  (1 + x)^{2(n-a_1- \dots -a_r-2-2^2- \dots 2^{t-1})}  \prod_{j=1}^{t-1} \left ( 1+ x^{2^j} \right )^2
\prod_{j=1}^r \left ( 1+ x^{a_j} \right )^2 
 }
{x^{n+1}} \right ] \quad .
$$
But (transferring a factor of $(1+x)^2$ from the second factor to the product,   $\prod_{j=1}^{t-1} \left ( 1+ x^{2^j} \right )^2$),
we have
$$
 (1 + x)^{2(n-a_1- \dots -a_r-2-2^2- \dots 2^{t-1})}  \prod_{j=1}^{t-1} \left ( 1+ x^{2^j} \right )^2
=
(1 + x)^{2(n-a_1- \dots -a_r-1-2-2^2- \dots 2^{t-1})}  \prod_{j=0}^{t-1} \left ( 1+ x^{2^j} \right )^2 \quad .
$$
Hence,
$$
A(\mu_0)(n)=
- \frac{1}{2} \, Coeff_{x^0} \left [
\frac{(1-x)^2  (1 + x)^{2(n-a_1- \dots -a_r-1-2-2^2- \dots 2^{t-1})}  \prod_{j=0}^{t-1} \left ( 1+ x^{2^j} \right )^2
\prod_{j=1}^r \left ( 1+ x^{a_j} \right )^2 
 }
{x^{n+1}} \right ] \quad .
$$
By Euler's good-old $(1-x)\prod_{j=0}^{t-1}(1+x^{2^j})=1-x^{2^t}$. Hence
$$
A(\mu_0)(n)=
-\, \frac{1}{2} \, Coeff_{x^0} \left [
\frac{(1-x^{2^t})^2  (1 + x)^{2(n-a_1- \dots -a_r-1-2-2^2- \dots 2^{t-1})}  \prod_{j=1}^r \left ( 1+ x^{a_j} \right )^2 }
{x^{n+1}} \right ] \quad .
$$
On the other hand, since $\mu'_0=Sort(a_1, \dots, a_r, 2^t)$, and all the $a_i$'s are odd, we have
$$
B(\mu'_0)(n+2)= - Coeff_{x^0}
\left [ 
\frac{(1+x)^{2n+2 -2(a_1 + \dots + a_r+2^t) }}{x^{n+1}} \cdot (x^{2^t}-1)^2 \cdot \prod_{j=1}^{r} (x^{a_j} +1)^2 
\right ] \quad .
$$
This completes the proof, since $-(1+2+ 2^2+  \dots +2^{t-1} ) = 1-2^t$ \quad . \halmos 

{\bf  Acknowledgment}

The research for this work was done while the second-named author visited the 
Faculty of Mathematics at the Weizmann Institute of Science,
during the week of Oct. 5-9, 2015. He wishes to thank the Weizmann Institute
for its hospitality, and its dedicated stuff, most notably Gizel Maimon.

{\bf References}

[M] I. G. Macdonald, {\it ``Symmetric Functions and Hall Polynomials''}, 2nd ed., Clarendon Press, Oxford, 1995.

[RRZ] Alon Regev, Amitai Regev, and Doron Zeilberger,
 {\it Identities in Character Tables of $S_n$}, J. Difference Equations and Applications, \hfill\break
DOI: 10.1080/10236198.2015.1081386, published online 11 Sep 2015, volume and page tbd.  \hfill\break
{\tt http://www.math.rutgers.edu/\~{}zeilberg/mamarim/mamarimhtml/sn.html} \quad .

[S] Neil Sloane, {\it The On-Line Encyclopedia of Integer Sequences}, {\tt http://oeis.org} \quad .

\bigskip
\hrule
\smallskip
Amitai Regev, Department of Pure Mathematics, Weizmann Institute of Science, Rehovot 76100, Israel ;
amitai dot regev at weizmann dot ac dot il \quad .
\smallskip
Doron Zeilberger, Department of Mathematics, Rutgers University (New Brunswick), Hill Center-Busch Campus, 110 Frelinghuysen
Rd., Piscataway, NJ 08854-8019, USA.  \hfill\break
zeilberg at math dot rutgers dot edu \quad .
\smallskip
\hrule
\medskip
Oct. 20, 2015.

\end